\documentclass[11pt]{article}      
\usepackage[centertags]{amsmath}
                      \usepackage{amsfonts,amssymb,hyperref,setspace}                                    

\newtheorem{theorem}{Theorem}
\newtheorem{remark}[theorem]{Remark}
\newtheorem{definition}[theorem]{Definition}
\newtheorem{corollary}[theorem]{Corollary}
\newtheorem{lemma}[theorem]{Lemma}
\newtheorem{assumption}[theorem]{Assumption}

\title{\LARGE \bf Further Results on Lyapunov Functions  for\\
Slowly Time-Varying Systems}

\author{ Fr\'{e}d\'{e}ric Mazenc\thanks{Mazenc is
with Projet MERE INRIA-INRA, UMR Analyse des Syst\`{e}mes et
Biom\'{e}trie INRA, 2, pl. Viala,
        34060 Montpellier, France (email:
        Frederic.Mazenc@ensam.inra.fr).} \and
Michael Malisoff\thanks{Malisoff is with the Department of
Mathematics, Louisiana State University, Baton Rouge, LA 70803-4918
USA (email: malisoff@lsu.edu).}}

\begin{document}

\maketitle

\begin{abstract}
We provide general methods for explicitly constructing strict
Lyapunov functions  for fully nonlinear slowly time-varying
 {}systems. Our results apply to cases where the given
dynamics and corresponding frozen dynamics are not necessarily
exponentially stable.  This complements  our previous Lyapunov
function constructions for rapidly time-varying dynamics.  We also
explicitly construct  input-to-state stable Lyapunov functions for
slowly time-varying control systems.
 We
illustrate our findings by constructing explicit Lyapunov functions
for a pendulum model, an example from identification theory, and a
perturbed friction model.  \medskip\medskip

\end{abstract}

\section{Introduction}
This paper is devoted to the study of fully nonlinear slowly
time-varying systems of the form
\begin{equation}
\label{stv} \dot x=f(x,t,t/\alpha)
\end{equation}
for large values of the constant $\alpha>0$ (but see Section
\ref{iss} below for the extension to systems with controls). See
Section \ref{defs} for our standing assumptions on (\ref{stv}). Such
systems arise in a large variety of important engineering
applications such as the control of friction and pendulums
\cite{K02, PA02b, S94}.  It is therefore of great interest in
control engineering to develop methods for determining whether
slowly time-varying systems   are uniformly globally asymptotically
stable (UGAS). When (\ref{stv}) is UGAS, it is also highly desirable
to have general methods for constructing explicit closed form
Lyapunov functions for (\ref{stv}).  See for example \cite{A99,
AS99,  ASW00, M03, MM06, MMdQ06} for discussions on the
essentialness of Lyapunov functions for feedback design and
robustness analysis.  See also \cite{MMdQ06,PA02a} for the dual
problem of stabilizing {\em rapidly} time-varying systems, and see
Remark \ref{transform} below for the relationship between our
methods for constructing Lyapunov functions for rapidly and slowly
time-varying systems.

One popular approach to studying (\ref{stv}) is to first establish
exponential stability of the corresponding ``frozen dynamics''
\begin{equation}
\label{frozen} \dot x=f(x,t,\tau)
\end{equation}
for all relevant values of the parameter $\tau$ including cases
where the exponent in the exponential decay estimate can be negative
or positive for some values of $\tau$ but is positive on average
\cite{K02, PA02b, S94}.  The stability of the frozen dynamics is
then used to establish stability of (\ref{stv}). However, these
earlier results do not lead to explicit Lyapunov functions for
(\ref{stv}) that would be needed for robustness analysis. The main
goals of our work are (i) to show that explicit Lyapunov functions
for (\ref{stv}) can be {}explicitly constructed in terms of a
suitable class of oftentimes  readily available Lyapunov functions
for (\ref{frozen}) {}when the constant $\alpha>0$ is large enough,
and (ii) to show how to relax the exponential like stability
assumptions on (\ref{frozen}) and also allow cases where $\tau$ is a
vector, thereby broadening the class of dynamics to which the frozen
dynamics method can be applied.

The rest of this paper is organized as follows.  In Section
\ref{defs} we provide the relevant definitions and standing
assumptions on (\ref{stv}).  We state and prove our main result in
Sections \ref{mains}-\ref{mainp}.  In Section \ref{nonex}, we extend
our main result to cases where the Lyapunov functions for the frozen
systems satisfy less restrictive properties than those {}in Sections
\ref{mains}-\ref{mainp}.
 In Section
\ref{ex}, we  illustrate {}the wide applicability of our results
using four examples. In the first two examples, the family of
Lyapunov like functions for (\ref{frozen}) are independent of
$\tau$, so the strict Lyapunov functions we construct for
(\ref{stv}) are valid for all $\alpha>0$ i.e. (\ref{stv}) is UGAS
for all $\alpha>0$. Our next two examples involve a mass spring
model with slowly time-varying coefficients from \cite{DDNZ00} and
an identification model similar to those studied in \cite{PA02a},
and illustrate the more general situation where
 (\ref{stv}) is not necessarily  UGAS for all values of
$\alpha>0$. In each case, the dynamics have slowly time-varying
coefficients and so are beyond the scope of the previously known
Lyapunov construction methods. In Section \ref{iss}, we show how to
extend our results to systems with controls using input-to-state
stability. We close in Section \ref{concl} with some suggestions for
further research.

\section{Definitions, Assumptions, and Lemmas}
\label{defs}

We let ${\cal K}_\infty$ denote the set of all continuous functions
$\rho:[0,\infty)\to[0,\infty)$ for which (i) $\rho(0)=0$ and  (ii)
$\rho$ is strictly increasing and unbounded.
We let ${\cal KL}$ denote the class of all continuous functions
$\beta:[0,\infty)\times [0,\infty)\to[0,\infty)$ for which
\begin{itemize}\addtolength{\itemsep}{-0.1\baselineskip}
\item[(I)]
 $\beta(\cdot, t)\in {\cal K}_\infty$ for each $t\ge 0$,
\item
[(II)] $\beta(s,\cdot)$ is non-increasing for each $s\ge 0$, and
\item[(III)] $\beta(s,t)\to 0$ as $t\to +\infty$ for each $s\ge 0$.
\end{itemize}When we say that a function $\rho$ is {\em smooth} (a.k.a.
$C^1$), we mean it is continuously differentiable, written $\rho\in
C^1$. (For functions $\rho$ defined on $[0,\infty)$, we interpret
$\rho'(0)$ as a one-sided derivative, and continuity of $\rho'$ at
$0$ as one-sided continuity.)  We let $|\cdot|$ denote the Euclidean
norm. A continuous function $\rho:[0,\infty)\to [0,\infty)$ is
called {\em positive definite} provided it is zero only at zero.
When $r\mapsto p(r)\in \mathbb{R}^d$ is a function with
differentiable components, we use $p'(r)$ to denote the vector
$(p_1'(r),....,p_d'(r))$.

The following definitions and lemma apply to general nonlinear
systems
\begin{equation}
\label{general} \dot x=h(x,t)
\end{equation}
evolving on the state space $\mathbb{R}^n$ where $h$ is locally
Lipschitz
 (but see Section
\ref{iss} for the extension to control systems). Later we specialize
to systems with multiple time scales and frozen parameters, e.g.
$h(x,t)=f(x,t,p(t/\alpha))$ or $h(x,t)=f(x,t,\tau)$ for given
constant parameters $\alpha$ and $\tau$ and suitable functions $p$.
We always assume (\ref{general}) is {\em forward complete} meaning
for each $x_o\in {\mathbb R}^n$ and $t_o\in {\mathbb R}_{\ge
0}:=[0,\infty)$ there exists a unique trajectory \[[t_0,\infty)\;
\ni\;  t\; \mapsto\;  \phi(t; t_o, x_o)\] for (\ref{general}) that
satisfies $x(t_o)=x_o$.  We assume all of our uncontrolled dynamics
(\ref{general}) are {\em uniformly state bounded}  meaning there
exists $\alpha_h\in \mathcal{K}_\infty$ such that $|h(x,t)|\le
\alpha_h(|x|)$ everywhere.

\begin{definition}
\label{ugasdef}  We say that (\ref{general}) is {\em uniformly
globally asymptotically stable (UGAS)} provided there exists
$\beta\in \mathcal{KL}$ such that {}\[\tag{UGAS}|\phi(t;t_o,x_o)|\le
\beta(|x_o|,t-t_o)\] for all  $x_o\in \mathbb{R}^n$, $t_o\in
[0,\infty)$, and  $t\ge t_o$.
\end{definition}
\begin{definition}
\label{lfdef} A smooth function $W:\mathbb{R}^n\times [0,\infty)\to
[0,\infty)$ is called a {\em Lyapunov function} for (\ref{general})
{}provided there are functions $\alpha_1,\alpha_2\in
\mathcal{K}_\infty$ and a positive definite function $\alpha_3$ such
that\smallskip
\begin{itemize}\item[]\begin{itemize}\addtolength{\itemsep}{.7\baselineskip}
\item[(L1)\ ]
$\alpha_1(|x|)\le W(x,t)\le \alpha_2(|x|)${}\ \ and
\item[(L2)\  ]
$W_t(x,t)+W_x(x,t)h(x,t)\le
-\alpha_3(|x|)$\end{itemize}\end{itemize}\smallskip hold for all
$t\ge 0$ and $x\in \mathbb{R}^n$.\end{definition}

The subscripts on $W$ denote partial gradients.  In what follows, we
often omit the arguments $x$, $t$, etc. in our functions when they
are clear from the context; and  all (in)equalities should be
interpreted to hold wherever they make sense.
 A smooth function $W:{\mathbb R}^n\times
[0,\infty)\to [0,\infty)$ that admits $\alpha_1,\alpha_2\in
\mathcal{K}_\infty$ such that (L1) holds everywhere is called {\em
uniformly proper and positive definite}. The following lemma is
standard \cite{ELW00, K02}:
\begin{lemma}\label{impl}
If (\ref{general}) admits a Lyapunov function, then it is UGAS.
\end{lemma}
A simple application of Fubini's Theorem yields the formula
\[\int_{t-c}^t\int_s^t\Theta(l)\, dl\, ds= \int_{t-c}^t
(r-t+c)\Theta(r)dr\] and therefore also the following \cite{MM06}:
\begin{lemma}\label{fubl}
Let $\Theta:\mathbb{R}\to \mathbb{R}$ be continuous and bounded in
norm by some constant $\bar M>0$,   and  $c>0$ be given.
Then\smallskip
\begin{itemize}\addtolength{\itemsep}{0.25\baselineskip}
\item[(A)]\ \
$\left\vert\displaystyle\int_{t-c}^{t}\int_{s}^{t} \Theta(l) dl
ds\right\vert\le \dfrac{c^2\bar M}{2}$\ \ \  and \item[(B)] \ \
$\dfrac{d}{dt}\displaystyle \int_{t-c}^{t}\int_{s}^{t}\Theta(l) dl
ds=c\Theta(t)-\int_{t-c}^t\Theta(r)dr$
\end{itemize}\smallskip
hold for all $t\in \mathbb{R}$.\end{lemma}

\section{Statement of Main Result and Remarks}
\label{mains}

For simplicity, we assume our system (\ref{stv}) has the form
\begin{equation}
\label{1} \dot{x} = f(x,t,p(t/\alpha))
\end{equation}
where $p:\mathbb{R}\to \mathbb{R}^d$ (for some integer $d$) is
bounded  and its components $p_1,\ldots, p_d$ have bounded first
derivatives.   We set
{}\[\bar p:=\sup\{|p'(r)|: r\in \mathbb{R}\}\; \; \; \; {\rm and}\;
\; \; \;  \mathcal{R}(p):=\{p(t): t\in \mathbb{R}\}.\] Our next
assumption is a variant of those  of \cite[Theorem 2]{PA02b} (but
see Section \ref{nonex} for results under weaker assumptions).

\begin{assumption}
\label{as1} There exist  $\alpha_1,\alpha_2\in \mathcal{K}_\infty$,
positive constants $c_a$, $c_b$, and $T$, a continuous function
$q:\mathbb{R}^d\to \mathbb{R}$,  and a $C^1$ function
$V:\mathbb{R}^n\times [0,\infty)\times \mathbb{R}^d\to [0,\infty)$
such that\smallskip
\begin{itemize}\addtolength{\itemsep}{0.25\baselineskip}
\item[$A_1$] $\alpha_1(|x|)\le V(x,t,\tau)\le \alpha_2(|x|)$,
\item[$A_2$] $V_t(x,t,\tau)+V_x(x,t,\tau)f(x,t,\tau)\le
-q(\tau)V(x,t,\tau)$,
\item[$A_3$] $|V_\tau(x,t,\tau)|\le c_aV(x,t,\tau)$, and
\item[$A_4$] $\int_{t-T}^tq(p(s))ds\ge c_b$
\end{itemize}\smallskip
hold for all $x\in\mathbb{R}^n$, $t\ge 0$, and $\tau\in {\mathcal
R}(p)$.
\end{assumption}

Note that $A_2$ is weaker than the standard exponential stability
property of the frozen dynamics since  we do not require $\alpha_1$
or $\alpha_2$ to be quadratic functions and moreover $q(\tau)$ can
take non-positive values for some choices of the vector parameter
$\tau$. However, $A_4$ requires that $q$ be positive on average
along the vector $p(s)$.

\begin{theorem}
\label{th1} If  (\ref{1}) satisfies Assumption \ref{as1}, then for
each constant $\alpha>2Tc_a\bar p/c_b$, the dynamics (\ref{1}) are
UGAS and
\begin{equation}
\label{5} \! \! \! V^\sharp_\alpha(t,x) :=
e^{\frac{\alpha}{T}\displaystyle\int_{\frac{t}{\alpha} -
T}^{\frac{t}{\alpha}} \displaystyle\int_{s}^{\frac{t}{\alpha}}
q(p(l)) dl\,  ds} V(x,t,p(t/\alpha))
\end{equation}
is a Lyapunov function for (\ref{1}).
\end{theorem}

\begin{remark}
\label{compare} Compared to the known results \cite{K02, PA02b,
S94}, the novelty of Theorem \ref{th1} is  that (a) we allow fully
nonlinear systems including cases where the function $q$ can take
both positive and negative values (which corresponds to the
allowance in \cite{S94} of eigenvalues that wander into the right
half plane while remaining in the strict left half plane on average)
and (b) we provide an explicit Lyapunov function (\ref{5}) for the
original slowly time-varying dynamics.  In general, the conclusion
of Theorem \ref{th1} may or may not hold for small values of
$\alpha$.  We illustrate {}this in Section \ref{ex} below.
\end{remark}

\begin{remark}
\label{transform} A completely different approach to  slowly
time-varying systems $\dot x=f(x,t,t/\alpha)$ (for large constants
$\alpha>0$) is to transform the system into a rapidly time-varying
system and to then try to construct a Lyapunov function for the
resulting rapidly time-varying system directly.  The transformation
is done by simply setting $s=t/\alpha$ which gives rise to the new
rapidly time-varying system
\begin{equation}
\label{newvar} \dot x(s)=g(x(s),s,\alpha s):=\alpha
f(x(s),s\alpha,s)\end{equation} in terms of the new rescaled time
variable $s$. However, it might be difficult to apply the Lyapunov
function construction methods of \cite{MMdQ06} or other known
methods to build an explicit Lyapunov function for (\ref{newvar}).
This is because these earlier results are for fast time-varying
dynamics having a different form from (\ref{newvar}) and moreover
they require Lyapunov functions for so-called limiting dynamics; cf.
\cite[Property 2]{PA02a}, and see \cite[Section 3.1]{MMdQ06} for a
generalization in the same vein.
 This motivates our direct
construction of Lyapunov functions for slowly time-varying dynamics,
which may be viewed as a complementary approach to the time
rescaling method since we do not require limiting dynamics.
\end{remark}

\section{Proof of Theorem \ref{th1}}
\label{mainp}  By $A_2$-$A_3$ and our choice of $\bar p$, the time
derivative  of
\begin{equation}
\label{v1def} \hat V(x,t) := V(x,t,p(t/\alpha))\end{equation} along
the trajectories of (\ref{1}) satisfies:
\[
\begin{array}{rcl}\renewcommand{\arraystretch}{1.5}
\! \! \! \! \! \!  \! \dot{\hat V} & = & \! \! \!
V_t(x,t,p(t/\alpha)) + V_x(x,t,p(t/\alpha)) f(x,t,p(t/\alpha))
+ V_\tau(x,t,p(t/\alpha))
\dfrac{p'(t/\alpha)}{\alpha}
\\
& \leq & \! \! \! - q(p(t/\alpha)) \hat V(t,x) +
V_\tau(x,t,p(t/\alpha))
 \dfrac{p'(t/\alpha)}{\alpha}\\ & \leq & \left[- q(p(t/\alpha))  +
\dfrac{c_a \bar p}{\alpha}\right] \hat V(x,t).
\end{array}\]
To simplify the notation, let us define
\begin{equation}\label{eta}
E(t,\alpha):= e^{\frac{\alpha}{T}\displaystyle\int_{\frac{t}{\alpha}
- T}^{\frac{t}{\alpha}}
\left[\displaystyle\int_{s}^{\frac{t}{\alpha}} q(p(l)) dl\right]
ds}.\end{equation} Since $p$ is bounded and $q$ is continuous,
\begin{equation}\label{theta}\Theta(t):=q(p(t))\end{equation} is bounded in norm by some constant $\bar
M>0$. Therefore,  along the trajectories of (\ref{1}), Lemma
\ref{fubl} (B) with the   choices (\ref{theta}) and $c=T$ gives
\[
\dot{V}^\sharp_\alpha =E(t,\alpha) \left[\dot{\hat V} +
\left\{q(p(t/\alpha)) -
\frac{1}{T}\displaystyle\int_{\frac{t}{\alpha} -
T}^{\frac{t}{\alpha}} q(p(l)) dl\right\} \hat V\right].
\]
Substituting the formula for $\dot {\hat V}$,  it follows from $A_4$
that
\begin{equation}\label{pre}\begin{array}{rcl} \dot{V}^\sharp_\alpha &\leq& E(t,\alpha)
\left[\dfrac{c_a \bar p}{\alpha}  -
\dfrac{1}{T}\int_{\frac{t}{\alpha} -
T}^{\frac{t}{\alpha}} q(p(l)) dl\right]\hat V\\[1em]
&\le & E(t,\alpha) \left[\dfrac{c_a\bar p}{\alpha} -
\dfrac{c_b}{T}\right] \hat V(x,t).
\end{array}\end{equation}
 Applying Lemma \ref{fubl} (A) with the   choices (\ref{theta}) and
 $c=T$ gives
\[e^{\alpha T\bar M/2}\ge E(t,\alpha)\ge e^{-\alpha T\bar M/2}\]
everywhere. Hence, for $\alpha>2Tc_a\bar p/c_b$, (\ref{pre}) gives
\begin{equation}\label{ly1}
\dot{V}^\sharp_\alpha(x,t)\; \; \le\; \;  -\frac{c_b}{2T}e^{-\alpha
T\bar M/2}\hat V(x,t)\le -\alpha_3(|x|),
\end{equation}
where $\alpha_3(s)=\frac{c_b}{2T}e^{-\alpha T\bar M/2}\alpha_1(s)$
is positive definite; and
\begin{equation}\label{ly2}
\begin{array}{rcl}
\hat \alpha_1(|x|)\; \; \le\; \;  V^\sharp_\alpha(x,t)\; \; \le\; \;
\hat \alpha_2(|x|)\end{array}
\end{equation}
everywhere, where \[\hat \alpha_1(s):=e^{-\alpha T\bar
M/2}\alpha_1(s),\; \; \; \hat \alpha_2(s):=e^{\alpha T\bar
M/2}\alpha_2(s)\] are of class $\mathcal{K}_\infty$.  Since
(\ref{ly1})-(\ref{ly2}) imply that $V^\sharp_\alpha$ is a Lyapunov
function for (\ref{1}), Theorem \ref{th1}  follows from Lemma
\ref{impl}.

\section{More General Families of Lyapunov Functions}
\label{nonex} We next show how to relax requirements $A_2$-$A_3$
from Theorem \ref{th1}.  We continue to use the notation we
introduced in Section \ref{mains}. We assume the following in the
rest of this section:

\begin{assumption}
\label{as2} There exist  $\tilde \alpha_1,\tilde \alpha_2\in
\mathcal{K}_\infty$, a positive definite $C^1$ function $\mu$,
positive constants $T$, $\tilde c_a$, and  $\tilde c_b$, a
continuous function $\tilde q:\mathbb{R}^d\to \mathbb{R}$, and a
$C^1$ function $\tilde V:\mathbb{R}^n\times [0,\infty)\times
\mathbb{R}^d\to [0,\infty)$ such that
\begin{equation}\label{mucon}
\displaystyle\lim_{r\to+\infty}\int_1^r\frac{1}{\mu(l)}dl=+\infty\end{equation}
and
\smallskip
\begin{itemize}\addtolength{\itemsep}{0.25\baselineskip}
\item[$\tilde A_1$] $\tilde \alpha_1(|x|)\le  \tilde V(x,t,\tau)\le \tilde
\alpha_2(|x|)$,
\item[$\tilde A_2$] $\tilde V_t(x,t,\tau)+\tilde V_x(x,t,\tau)f(x,t,\tau)\le
-\tilde q(\tau)\mu(\tilde V(x,t,\tau))$,
\item[$\tilde A_3$] $|\tilde  V_\tau(x,t,\tau)|\le \tilde c_a\mu(\tilde
V(x,t,\tau))$, and
\item[$\tilde A_4$] $\int_{t-T}^t\tilde q(p(s))ds\ge \tilde c_b$
\end{itemize}\smallskip
hold for all $x\in\mathbb{R}^n$, $t\ge 0$, and $\tau\in {\mathcal
R}(p)$.
\end{assumption}

Notice that  Assumption \ref{as1} is the special case of Assumption
\ref{as2} in which  $\mu(l)\equiv l$. We prove the following:

\begin{theorem}
\label{th2} If  (\ref{1}) satisfies Assumption \ref{as2}, then there
exists $k\in \mathcal{K}_\infty$ for which the requirements of
Assumption \ref{as1} are satisfied with $V:=k(\tilde V)$. Therefore,
for each sufficiently large value of the constant $\alpha>0$, the
dynamics (\ref{1}) is UGAS and admits a Lyapunov function of the
form (\ref{5}).
\end{theorem}

It suffices to prove the first statement of Theorem \ref{th2} since
the second statement is immediate  from Theorem \ref{th1}.  To this
end, we use the following {}important observation:

\begin{lemma}
\label{ll} If $\mu\in C^1$ is positive definite,  then
\begin{equation}\label{newkey}\lim_{r\to
0^+}\int_1^r\frac{1}{\mu(l)}dl=-\infty.\end{equation}
\end{lemma}\smallskip
\noindent{\em Proof:}\ \  Since $\mu\in C^1$ and $\mu(0)=0$, we can
find $c_3>0$ such that $\mu(r)\le c_3r$ for all $r\in [0,1]$. Hence,
for each $r\in (0,1]$, we get
\begin{equation}
\label{cz13}
\begin{array}{rcl}
\displaystyle\int_{r}^{1}\frac{1}{\mu(l)} dl & \geq &
\displaystyle\int_{r}^{1}\frac{1}{c_3 l} dl = - \dfrac{1}{c_3}
\ln(r)
\end{array}
\end{equation}
It follows that, for all $r \in (0,1]$,
\begin{equation}
\label{cz14}
\begin{array}{rcl}
\displaystyle\int_{1}^{r}\frac{1}{\mu(l)} dl & \leq & \dfrac{1}{c_3}
\ln(r).
\end{array}
\end{equation}
Since $\lim_{r \rightarrow 0^+} \ln(r) = - \infty$, the lemma
follows. \hfill$\blacksquare$

Given a constant $\xi>0$ which we specify later, it follows from
(\ref{mucon}) and  {(\ref{newkey}) that
\begin{equation}
\label{cz15}
\begin{array}{rcl}
k(r) & = & \left\{\begin{array}{lcl}
e^{\xi\displaystyle\int_{1}^{r}\frac{1}{\mu(l)} dl},&& r>0\\
0,&& r=0\end{array}\right.
\end{array}
\end{equation}
 is
continuous and unbounded.  {}In particular, $k(r)\to 0$ as $r\to
0^+$. Set $B = \sup\{\mu'(s):0\le s\le 1\}$, which is positive since
$\mu$ is positive definite.
\begin{lemma}\label{cone}  The function (\ref{cz15}) is $C^1$ when  $\xi=2B$.\end{lemma}\noindent{\em Proof:}\ \  It
suffices to prove that \begin{equation} \label{goal} k'(r)\to 0 \;
\; {\rm as}\; \;  0<r\to 0^+,\end{equation} since $k(0)=0$, because
then $k(r)/r\to 0=k'(0)$ as $r\to 0^+$.
 To this end, first note that
\begin{equation}
\label{sz1} k'(r)  = \frac{\xi}{\mu(r)} e^{-
\xi\displaystyle\int_{r}^{1}\frac{1}{\mu(l)} dl}, \; \; \forall r>0
\end{equation}
 and  that for all $r\in (0,1]$, we have
\[
\frac{1}{B} \left[\ln(\mu(1)) - \ln(\mu(r))\right] \; =\;
\displaystyle\int_{r}^{1} \frac{\mu'(l)}{B \mu(l)} dl\; \le\;  
\displaystyle\int_{r}^{1}\frac{1}{\mu(l)} dl
\]
by our choice of $B$.  Since $\mu$ is positive definite and $\xi$ is
positive, this implies \[ \frac{\xi}{\mu(r)} e^{\frac{- \xi}{B}
\left[\ln(\mu(1)) - \ln(\mu(r))\right]} \; \geq\; \frac{\xi}{\mu(r)}
e^{- \xi \displaystyle\int_{r}^{1}\frac{1}{\mu(l)} dl}
\; =\; k'(r)
\] for all $r\in
(0,1]$, i.e.,
\begin{equation}
\label{ez8} \xi (\mu(1))^{\frac{- \xi}{B}} \mu(r)^{\frac{\xi}{B}-1}
\; \; \geq\; \;  k'(r) {}\; \; \ge\;  \; 0\;  \; \; \; \forall r \in
(0,1]
\end{equation}
so (\ref{goal}) follows from our choice of $\xi$.
\hfill$\blacksquare$

The fact that $k\in \mathcal{K}_\infty\cap C^1$ is now immediate
from Assumption \ref{as2} and Lemmas \ref{ll} and \ref{cone}.  Let
us now verify that $k$ satisfies the requirements of Theorem
\ref{th2}. From the definition of $k$, we deduce that
\[
k'(\tilde V)\mu(\tilde V) \; =\;  \frac{2 B}{\mu(\tilde V)} e^{- 2
B\displaystyle\int_{\tilde V}^{1}\frac{1}{\mu(l)} dl} \mu(\tilde
V)\; =\;  2 B k(\tilde V)\] when $\tilde V\ne 0$. Therefore, by
assumption $\tilde A_2$, the time derivative of $V=k(\tilde V)$
along the trajectories of (\ref{1}) satisfies
\[\renewcommand{\arraystretch}{1.5}\begin{array}{rcl}
\! \! \! V_t(x,t,\tau)+V_x(x,t,\tau)f(x,t,\tau)
&=& k'(\tilde V(x,t,\tau))[\tilde V_t(x,t,\tau)+\tilde
V_x(x,t,\tau)f(x,t,\tau)]\\
&\le& - k'(\tilde V(x,t,\tau))\tilde q(\tau)\mu(\tilde
V(x,t,\tau))\\
&=&- 2 B \tilde q(\tau) V(x,t,\tau)\end{array}\] everywhere, and
condition $\tilde A_3$ from Assumption \ref{as2} implies
\[|V_\tau|\; =\;  k'(\tilde V) |\tilde V_\tau| \; \leq\;  \tilde c_a k'(\tilde
V) \mu(\tilde V)  \; =\;  2B\tilde c_a V\] everywhere. Therefore
Assumption \ref{as1} holds using
\begin{equation}
\label{choicess} \alpha_i(s):=k\circ \alpha_i(s)\; \;  {\rm for}\;
\;  i=1,2, \; \; \; \; V:=k(\tilde V),\; \; \;  \; q(\tau) := 2 B
\tilde q(\tau),\; \; \;  \; c_a:=2B\tilde c_a, \; \; \; \;
c_b:=2B\tilde c_b.\; \; \end{equation} The result now follows from
Theorem \ref{th1}.

\section{Examples}
We illustrate our constructions using four examples.  In the first
two  examples, the functions $V$ from Assumption \ref{as1} do not
depend on the frozen parameter $\tau$, so we can conclude that
(\ref{1}) is UGAS for all values of the constant $\alpha>0$.  We
then turn to a slowly time-varying friction dynamics  and an example
from identification where $V$ depends on $\tau$, and where we can
consequently only conclude the UGAS property of (\ref{1}) when
$\alpha>0$ is sufficiently large. Set $\dot
V(x,t,\tau):=V_t(x,t,\tau)+V_x(x,t,\tau)f(x,t,\tau)$ everywhere.

\label{ex} \subsection{Stability for all  {}$\mathbf{\alpha>0}$: a
scalar example}

Consider the one-dimensional system
\begin{equation}
\label{16} \dot{x} = \frac{x}{\sqrt{1 + x^2}} \left[1 - 90
\cos^2\left(\frac{t}{\alpha}\right)\right]
\end{equation}
and the uniformly proper and positive definite function
\begin{equation}
\label{17} V(x,t,\tau) \equiv \bar V(x):=e^{\sqrt{1 + x^2}} - e.
\end{equation}
Let us verify Assumption \ref{as1} for this choice of $V$ and the
frozen dynamics
\begin{equation}
\label{21} \dot{x} = f(x,t,\tau):=\frac{x}{\sqrt{1 + x^2}} \left[1 -
90 \tau\right],\; \; \; 0\le \tau\le 1.
\end{equation}
This gives
\begin{equation}
\label{18}
\begin{array}{rcl}
\dot V(x,t,\tau)
& = & e^{\sqrt{1 + x^2}} \dfrac{x^2}{1 + x^2}\left[1 - 90
\tau\right]
\\[1em]
& = & e^{\sqrt{1 + x^2}} \dfrac{x^2}{1 + x^2} - 90 \tau e^{\sqrt{1 +
x^2}} \frac{x^2}{1 + x^2}.
\end{array}
\end{equation}
Simple calculus calculations everywhere give
\begin{equation}
\label{19} \frac{2 e^{\sqrt{2}}}{e - 1}\bar V(x) \; \geq\;
\frac{x^2}{1 + x^2} e^{\sqrt{1 + x^2}} \; \geq\;  \frac{1}{2}\bar
V(x)
\end{equation}
so (\ref{18}) everywhere gives
\begin{equation}
\label{20}
\begin{array}{rcl}
\dot{V}(x,t,\tau) & \leq & \left[\dfrac{2 e^{\sqrt{2}}}{e - 1} -
45\tau\right] \bar V(x).
\end{array}
\end{equation}
Moreover, for each $t\ge 0$, we get \[\int_{t-\pi}^{t} \left[45
\cos^2(s) - \frac{2 e^{\sqrt{2}}}{e - 1}\right] ds \; =\;
\pi\left(\frac{45}{2} - \frac{2 e^{\sqrt{2}}}{e - 1}\right) > 0
\]
which shows that Assumption \ref{as1} is satisfied.  We conclude
from Theorem \ref{th1}  that for large enough constants $\alpha>0$,
(\ref{16}) is UGAS and has the Lyapunov function
\begin{equation}
\label{23} \begin{array}{l}
e^{\frac{\alpha}{\pi}\displaystyle\int_{\frac{t}{\alpha} -
\pi}^{\frac{t}{\alpha}}\left[
\displaystyle\int_{s}^{\frac{t}{\alpha}} \left[45 \cos^2(l) -
\frac{2 e^{\sqrt{2}}}{e - 1}\right] dl \right] ds} \! \! \! \bar
V(x)
 \; \; =\; \;   e^{45\frac{\alpha}{4}\left[\sin(\frac{2 t}{\alpha}) + \pi -
\frac{4 \pi e^{\sqrt{2}}}{45(e - 1)}\right]}[e^{\sqrt{1 + x^2}} - e]
\end{array}
\end{equation}
where $\bar V$ is in (\ref{17}). In fact, since in this case $V$
does not depend on $\tau$, it follows from our proof of Theorem
\ref{th1} that for any constant $\alpha
> 0$, the system (\ref{16}) is UGAS and admits the Lyapunov function
(\ref{23}). \begin{remark} \label{91} The dynamics (\ref{16})
illustrates the fact that our approach applies to systems which are
not globally exponentially stable. Indeed, it is clear that
(\ref{16})  is not  globally exponentially stable since its vector
field is bounded in norm by the constant $91$.\end{remark}

\subsection{Stability for all {}$\mathbf{\alpha>0}$:  a pendulum example}

Our constructions also apply to the slowly time-varying pendulum
dynamics \cite{PA02b}
\begin{equation}
\label{pend}
\begin{array}{l}
\dot x_1= x_2{}\\[.2em]\dot x_2= -x_1-[1+b_2(t/\alpha)m(x,t)]x_2\end{array}
\end{equation}
assuming
\begin{itemize}
\item[$(\mathcal{P}1)$]$m:\mathbb{R}^2\times \mathbb{R}\to [0,1]$ is
Lipschitz  continuous; and {}\item[$(\mathcal{P}2)$]
$b_2:\mathbb{R}\to (-\infty,0]$ is globally bounded, and there are
positive  constants $T$ and $c_b$ such that $5+T\int_{t-T}^tb_2(l)\,
dl\ge c_b$ for all $t\in \mathbb{R}$.
\end{itemize}
The dynamics (\ref{pend}) was shown to be UGAS for certain choices
of the function $b_2$ in \cite{PA02a};  see \cite{K02} for related
results, and \cite{S94} for results that are restricted to the
linear case.  However, these earlier results do not lead to explicit
Lyapunov functions for (\ref{pend}). In order to build Lyapunov
functions for (\ref{pend}) for large constants $\alpha>0$, we use
the following observation for the corresponding frozen dynamics
$f(x,t,\tau):=(x_2,-x_1-[1+\tau m(x,t)]x_2)$:
\begin{lemma}\label{penlm}
The function $V(x):=x^2_1+x^2_2+x_1x_2$ satisfies $\nabla
V(x)f(x,t,\tau)\le -[1+5\tau]V(x)$ for all $x\in \mathbb{R}^n$,
$t\in \mathbb{R}$, and $\tau\le 0$.
\end{lemma}
\noindent{\em Proof:}\ By grouping terms, one easily shows that
\begin{equation}
\label{easily} \nabla V(x)f(x,t,\tau)=-V(x)-2\tau m(x,t)x^2_2-\tau
m(x,t)x_1x_2\end{equation} everywhere.  Since \[V\; \ge\;
x^2_1+x^2_2-|x_1x_2|\; \ge\; \dfrac{1}{2}x^2_1+\frac{1}{2}x^2_2\;
\ge\;  |x_1x_2|\] everywhere, we get
\begin{equation}
\label{putin} -2\tau m(x,t) x^2_2\le -4\tau m(x,t)V(x),\; \; \; \;
-\tau m(x,t)x_1 x_2\le -\tau m(x,t)V(x)\end{equation} everywhere.
The lemma follows by substituting (\ref{putin}) into (\ref{easily})
and recalling that $0\le m(x,t)\le 1$ everywhere.
\hfill$\blacksquare$

The following is an immediate consequence of Lemma \ref{penlm}, the
proof of Theorem \ref{th1}, and the fact that
$V(x):=x^2_1+x^2_2+x_1x_2$ only depends on $x$:

\begin{theorem}
\label{penthm} Let (\ref{pend}) satisfy
$(\mathcal{P}1)$-$(\mathcal{P}2)$.  Then (\ref{pend}) has the
Lyapunov function
\begin{equation}
\label{penlf} V^\sharp_\alpha(t,x) :=
e^{\frac{5\alpha}{T}\displaystyle\int_{\frac{t}{\alpha} -
T}^{\frac{t}{\alpha}} \displaystyle\int_{s}^{\frac{t}{\alpha}}
b_2(l)dl\, ds} (x^2_1+x^2_2+x_1x_2)
\end{equation}
for each choice of the constant $\alpha>0$.  Hence, (\ref{pend}) is
UGAS for all constants $\alpha>0$.
\end{theorem}

\subsection{Friction example revisited}
\label{revisit} We next illustrate Theorem \ref{th1} using the one
degree-of-freedom mass-spring system \cite{DDNZ00, MMdQ06}.  The
mass-spring system arises in the control of mechanical systems in
the presence of friction. However, in contrast to \cite{MMdQ06}
where the dynamics are assumed to be rapidly time-varying, here we
consider the case where the dynamics are {\em slowly} time-varying.
While slowly time-varying dynamics can be transformed into rapidly
time-varying dynamics by rescaling time, doing so for the slowly
time-varying mass spring system produces a new dynamic that does not
lend itself to the known methods; see Remark \ref{transform} above
for details. For this reason, we directly apply the slowly
time-varying theory we developed in the preceding sections.

Let us recall the model \cite{MMdQ06}.  The dynamics are given by
\begin{equation}
\label{frictionexample}
\begin{array}{rcl}
\! \! \! \dot{x}_{1} &=&x_{2}\\ \dot{x}_{2} &=&-\sigma
_{1}(t/\alpha)x_{2}-k(t)x_{1} -\left\{ \sigma _{2}(t/\alpha)+\sigma
_{3}(t/\alpha)e^{-\beta _{1}\mu (x_{2})}\right\} {\rm sat}
(x_{2})\end{array}\end{equation} where $x_{1}$ and $x_{2}$ are the
mass position and velocity, respectively; $\sigma
_{i}:[0,\infty)\to(0,1]$, $i=1,2,3$ denote positive time-varying
viscous, Coulomb, and static friction-related coefficients,
respectively; $\beta _{1}$ is a positive constant corresponding to
the Stribeck effect; $\mu (\cdot )$ is a positive definite function
also related to the Stribeck effect; $k$ denotes a positive
time-varying spring stiffness-related coefficient; and sat$(\cdot )$
denotes any continuous function having these properties:
\begin{equation}\label{sats}\renewcommand{\arraystretch}{1.2}
\begin{array}{l} \! \! \! \! \!
\! {\rm (a)}\; \; {\rm sat}(0)=0, \ \ \  {\rm (b)}\; \;  \xi \, {\rm
sat}(\xi )\geq 0\; \; \forall \xi\in {\mathbb R},\\ \! \! \! \! \!
\! {\rm (c)}\; \; \lim\limits_{\xi \rightarrow + \infty }{\rm
sat}(\xi )=+ 1,  \; \; \; {\rm (d)}\; \; \lim\limits_{\xi
\rightarrow - \infty }{\rm sat}(\xi )=- 1
\end{array}\end{equation} Following \cite{MMdQ06},
we model the saturation differentiably as
\begin{equation}
{\rm sat}(x_{2})=\tanh (\beta _{2}x_{2}),  \label{sat}
\end{equation}
where $\beta _{2}$ is a large positive constant, so $|{\rm
sat}(x_{2})|\le \beta_2|x_2|$ for all $x_2\in {\mathbb R}$. However,
unlike \cite{MMdQ06}, we assume the friction coefficients vary in
time {\em slower} than the spring stiffness coefficient so we
restrict to cases where $\alpha>1$.  We are going to establish the
stability of (\ref{frictionexample}) and construct corresponding
Lyapunov functions $V_\alpha$ when the constant $\alpha>0$ is
sufficiently large.

Our precise mathematical assumptions on  (\ref{frictionexample})
are: $k$ and the $\sigma_i$'s are  $C^1$ functions with bounded
derivatives; $\mu $ has a globally bounded derivative; and there
exist constants $c_b,T>0$ such that
\begin{equation}
\label{key1} \int_{t-T}^t\sigma_1(r)dr\; \ge\;  c_b \; \; \; \;
\forall t\ge 0.\end{equation}  We also assume this additional
condition whose physical interpretation is that the spring stiffness
is nonincreasing:
\[
\exists k_o, \bar k>0\; \; {\rm s.t.}\; \; k_o \le k(t)\le \bar k\;
\; {\rm and}\; \; k'(t)\le 0 \; \; \forall t\ge 0.
\]

The frozen dynamics $\dot x=f(x,t,\tau)$ for (\ref{frictionexample})
are

\begin{equation} \label{frictionexamplea}
\begin{array}{rcl}
\dot{x}_{1} &=&x_{2} \\
\dot{x}_{2} &=&-\tau_1 x_{2}-k(t)x_{1}-\left\{
\tau_2+\tau_3e^{-\beta _{1}\mu (x_{2})}\right\} {\rm sat}
(x_{2})\end{array}
\end{equation}
where $\tau=(\tau_1,\tau_2,\tau_3)\in [0,1]^3$ is now a {\em vector}
of parameters.  We apply our construction from Theorem \ref{th1}
with $p(t)= (\sigma_1(t), \sigma_2(t),\sigma_3(t))$ and the function
\begin{equation}
\label{Vchoice}
\begin{array}{l}
V(x,t,\tau)=A(k(t)x^2_1+x^2_2)+\tau_1x_1x_2\; \; \;  {\rm where}\;
\; \;  A=1+\dfrac{k_o}{2}+\dfrac{(1+2\beta_2)^2}{k_o}.\end{array}
\end{equation}

We first verify the conditions of Assumption \ref{as1}.  Since $A\ge
\max\{1, 1/k_o\}$ and $\tau_1\le 1$, we have {}\begin{equation}
\label{newbounds} \frac{1}{2}(x^2_1+x^2_2)\, \le\,  V(x,t,\tau)\,
\le\, A^2\bar k(|x_1|+|x_2|)^2\, \le\,  2A^2\bar k|x|^2
\end{equation}
everywhere.  Let us now compute $\dot V(x,t,\tau)$ for all values
$\tau\in [0,1]^3$. Since $k'(t)\le 0$ everywhere, this gives
\[\renewcommand{\arraystretch}{1.25}
\begin{array}{rcl}
\dot V(x,t,\tau) &\le&  V_x(x,t,\tau)f(x,t,\tau)
\\&=& [2Ak(t)x_1+\tau_1x_2]x_2-[2Ax_2+\tau_1x_1]\{\tau_1x_{2}+ \left[\tau_2+\tau _{3}e^{-\beta _{1}\mu
(x_{2})}\right] {\rm sat} (x_{2})+k(t)x_{1}\}.
\end{array}\]
Therefore, by grouping and canceling terms, we also have
\[
\renewcommand{\arraystretch}{2}
\begin{array}{rcl}
\dot V(x,t,\tau)&\le&
-\tau_1k_0x^2_1-(2A\tau_1-\tau_1)x^2_2+\tau_1(1+2\beta_2)|x_1x_2|\\
&\le&
-\tau_1\dfrac{k_o}{2}|x|^2-\left[\tau_1\dfrac{k_o}{2}x^2_1+(A-1/2)\tau_1x^2_2-
\tau_1(1+2\beta_2)|x_1x_2|\right]\\
&=&
-\tau_1\dfrac{k_o}{2}|x|^2-\tau_1\dfrac{k_o}{2}\left(|x_1|-\dfrac{1+2\beta_2}{k_o}|x_2|\right)^2\!
+\left(\dfrac{\tau_1(1+2\beta_2)^2}{2k_o}+\dfrac{\tau_1}{2}
-A\tau_1\right)x^2_2\\
&\le& -\dfrac{\tau_1k_o}{4A^2\bar k} V(x,t,\tau)
\end{array}\]
where the first inequality follows from (\ref{sats})(b), the
inequality $|{\rm sat}(x_2)|\le \beta_2|x_2|$, and the fact that
$\tau_i\in [0,1]$ for each $i$; the second inequality uses the fact
that $A-1/2\ge k_o/2$; and the last inequality is from the choice of
$A$ {}and the bounds (\ref{newbounds}).   Hence, Assumption
\ref{as1} of Theorem \ref{th1} readily follows from (\ref{key1})
with the choices
\[q(\tau)=\dfrac{\tau_1k_o}{4A^2\bar k},\; \; \; \; c_a=1.\] We conclude as
follows:

\begin{corollary}
\label{cor1} Under the preceding assumptions, there exists a
constant $\alpha_o>0$ such that for all {}constants
$\alpha>\alpha_o$, the system (\ref{frictionexample}) is UGAS and
admits the Lyapunov function
\begin{equation}
\label{fl}
\begin{array}{l}
\! \! \! \! \! \! \! \! \! V_\alpha(t,x) := V(x,t,p(t/\alpha))\,
e^{\frac{\alpha\bar b}{T}\displaystyle\int_{\frac{t}{\alpha} -
T}^{\frac{t}{\alpha}}
\displaystyle\int_{s}^{\frac{t}{\alpha}}\sigma_1(l) dl\,  ds}
\end{array}\end{equation} where $V$ is the function defined in  (\ref{Vchoice}),
$\bar b=k_o/(4A^2\bar k)$,
 and
$p(t)=(\sigma_1(t), \sigma_2(t),\sigma_3(t))$.
\end{corollary}

\subsection{Identification dynamics revisited}
Our Lyapunov function constructions also apply to the slowly
time-varying dynamics
\begin{equation}
\label{slowid} \dot x=h(t/\alpha)m(t)m^\top(t)x, \; \; x\in
\mathbb{R}^n
\end{equation}
assuming there are positive constants $T$, $\tilde c$, $\underline
\alpha$, and $\bar \alpha$ such that
\begin{itemize}\addtolength{\itemsep}{-0.25\baselineskip}
\item[$(\mathcal{I}1)$]
$h:\mathbb{R}\to[-\bar\alpha,0]$ is continuous with a  bounded first
derivative and $\int_{t-T}^th(r)dr\le -\underline \alpha$ for all
$t\in \mathbb{R}$.
\item[$(\mathcal{I}2)$]
$m:\mathbb{R}\to\mathbb{R}^n$ is continuous and satisfies
$|m(t)|\equiv 1$ and $\underline \alpha I\le \int_{t}^{t+\tilde c}
m(r)m^\top(r)dr\le \bar \alpha I$ for all $t\in \mathbb{R}$.
\end{itemize}
where $I$ is the identity matrix and for matrices $A,B\in
\mathbb{R}^{n\times n}$ we use $B-A\ge 0$ to mean that $B-A$ is
positive semi-definite.  {}See also  Remark  \ref{lastrm} below  for
the generalization of our result to control affine systems \[\dot
x\; \; =\; \; h(t/\alpha)m(t)m^\top(t)x\; +\; g(x,t,t/\alpha)u\] for
suitable matrix valued functions $g$.  The particular case $\dot
x=-m(t)m^\top(t)x$ of (\ref{slowid}) has been studied extensively in
the context of identification theory \cite{MMdQ06, PA02a}. In
\cite{MMdQ06}, we showed how to construct explicit Lyapunov
functions for the rapidly time-varying system $\dot x=f(\alpha
t)m(t)m^\top(t)x$ for appropriate nonpositive functions $f$ and
large positive constants $\alpha$.  However, these earlier results
do not lead to explicit Lyapunov functions for the slowly
time-varying dynamics (\ref{slowid}) for large constants $\alpha>0$;
see Remark \ref{transform} above. Instead, we construct Lyapunov
functions for (\ref{slowid}) using the following analogue of
\cite[Lemma 6]{MMdQ06}:
\begin{lemma}
\label{idlm} Assume there are positive constants $T$, $\tilde c$,
$\underline \alpha$, and $\bar \alpha$ such that
$(\mathcal{I}1)$-$(\mathcal{I}2)$ are satisfied and set
\begin{equation}
\label{newpchoice}
\begin{array}{l}
P(t,\tau)=\kappa I-\tau \displaystyle\int_{t-\tilde c}^t\int_s^t
m(l)m^T(l)\, {\rm d}l\, {\rm d}s,\; \; {\rm where}\; \;
\kappa=\dfrac{\tilde c}{2}+\dfrac{\bar \alpha^2\tilde
c^4}{2\underline \alpha}+{}\tilde c^2.
\end{array}
\end{equation}
Then for each $\tau\in [-\bar \alpha,0]$ the function
\[ V(x,t,\tau)\; =\; x^\top P(t,\tau)x\]  satisfies the requirements of
Assumption \ref{as1} for the frozen dynamics $f(x,t,\tau)=\tau
m(t)m^\top (t)x$ and $p(s)=h(s)$.\end{lemma}

\noindent{\em Proof:}\  We first apply Lemma \ref{fubl}(B) to
$\Theta(t)=x^\top m(t)m^\top(t)x$ for each $x$ and $\tau\in [-\bar
\alpha,0]$ to get
\[
\frac{\partial V}{\partial t}(x,t,\tau)=-\tau\tilde c x^\top
m(t)m^\top (t)x+\tau x^\top \left[\int_{t-\tilde c}^t m(l)m^\top
(l)dl\right]x\]
 and
 \[
\frac{\partial V}{\partial x}(x,t,\tau)f(x,t,\tau)=2\tau
x^\top\left[\kappa I-\tau\int_{t-\tilde c}^t\int_s^t m(l)m^T(l)\,
{\rm d}l\, {\rm d}s\right]m(t)m^\top (t)x
\]
so the derivative $\dot V=\frac{\partial V}{\partial
t}+\frac{\partial V}{\partial t} f$ along the trajectories of $f$
satisfies
\begin{equation}\label{dvt}
\begin{array}{rcl}
\dot V&\le & \tau[(2\kappa-\tilde c)|m^\top (t)x|^2+\underline
\alpha |x|^2]+ \tau^2|x||m^\top (t)x|\tilde c^2\\[.5em]
&=& \tau\{(2\kappa-\tilde c)|m^\top (t)x|^2+\underline \alpha |x|^2+
\tau|x||m^\top (t)x|\tilde c^2\}\end{array} \end{equation} where the
inequality follows from Lemma \ref{fubl}(A), $(\mathcal{I}2)$, and
the facts that $|m(t)|\equiv 1$ and $\tau\le 0$. Set
\[
\omega=\frac{\underline \alpha}{2\tilde c^2\bar \alpha}.\] Then the
triangle inequality gives
\[
|m^\top(t)x||x|\; \le\; \omega|x|^2+\frac{1}{4\omega}|m^\top(t)x|^2
\]
so since $\tau\le 0$, we get
\[\renewcommand{\arraystretch}{1.5}
\begin{array}{l}
(2\kappa-\tilde c)|m^\top(t)x|^2+\underline \alpha
|x|^2+\tau|x||m^\top(t)x|\tilde
c^2\\
\ge (2\kappa-\tilde c)|m^\top(t)x|^2+\underline \alpha
|x|^2+\omega\tau\tilde c^2|x|^2+\dfrac{\tau\tilde
c^2}{4\omega}|m^\top(t)x|^2\\
= \left[2\kappa-\tilde c+\dfrac{\tau\tilde
c^2}{4\omega}\right]|m^\top(t)x|^2+(\underline
\alpha+\omega\tau\tilde c^2)|x|^2\\
\ge \dfrac{\underline \alpha|x|^2}{2}
\end{array}
\]
by our choices of $\kappa$ and $\omega$.  This and (\ref{dvt}) give
$\dot V\le \frac{\tau\underline \alpha}{2}|x|^2$ everywhere; and
Lemma \ref{fubl}(A), the fact that $|m(t)|\equiv 1$,  and our choice
of $\kappa$ give
\[
\begin{array}{l}
\kappa|x|^2\; \le\;  V(x,t,\tau)\; \le\;  (\kappa +\tilde c^2\bar
\alpha)|x|^2\\[.2em]
\left\vert\dfrac{\partial V}{\partial\tau}(x,t,\tau)\right\vert\;
\le\;  \dfrac{{}\tilde c^2}{2}|x|^2\le V(x,t,\tau).
\end{array}\]

Hence,
\[
\dot V(x,t,\tau)\; \le\; -q(\tau)V(x,t,\tau), \; \; {\rm where}\; \;
q(\tau):=-\frac{\tau\underline\alpha}{2(\kappa+\tilde c^2\bar
\alpha)}
\]
everywhere. Therefore, we can satisfy the requirements of Assumption
\ref{as1} using $p(s)=h(s)$ and
\begin{equation}\label{endchoice}
\alpha_1(s)=\kappa s^2,\; \; \;  \alpha_2(s)=(\kappa +\tilde c^2\bar
\alpha)s^2,\; \; c_a=1,\; \; c_b=\frac{\underline
\alpha^2}{2(\kappa+\tilde c^2\bar \alpha)}.
\end{equation}
 This proves the
lemma.\hfill$\blacksquare$

The following is an immediate consequence of the preceding lemma and
Theorem \ref{th1}:

\begin{theorem}
\label{idthm} Let (\ref{slowid}) admit  positive constants $T$,
$\tilde c$, $\underline \alpha$, and $\bar \alpha$ such that
$(\mathcal{I}1)$-$(\mathcal{I}2)$ are satisfied and choose $c_b$ as
{}in (\ref{endchoice}). Then for any constant $\alpha \ge
2T\sup\{|h'(r)|: r\in \mathbb{R}\}/c_b$, the function
\begin{equation}
\label{idlf} V^\sharp_\alpha(t,x) := e^{-\frac{\alpha\underline
\alpha}{2T(\kappa+\tilde c^2\bar
\alpha)}\displaystyle\int_{\frac{t}{\alpha} - T}^{\frac{t}{\alpha}}
\displaystyle\int_{s}^{\frac{t}{\alpha}} h(l)dl\, ds} V(x,t,
h(t/\alpha))
\end{equation}
with $V(x,t,\tau)=x^\top P(t,\tau)x$ as defined in Lemma \ref{idlm}
is a Lyapunov function for (\ref{slowid}).  Hence, (\ref{slowid}) is
UGAS for all constants $\alpha \ge 2T\sup\{|h'(r)|: r\in
\mathbb{R}\}/c_b$.
\end{theorem}

\section{Input-to-State Stability}
\label{iss} We next extend our results to control affine slowly
time-varying systems
\begin{equation}
\label{controlled} \dot{x} = f(x,t,p(t/\alpha))+g(x,t,p(t/\alpha))u
\end{equation}
evolving on $\mathbb{R}^n$ with control values $u\in\mathbb{R}^m$,
where $f:\mathbb{R}^n\times [0, \infty)\times \mathbb{R}^d\to
\mathbb{R}^n$ and $g:\mathbb{R}^n\times [0, \infty)\times
\mathbb{R}^d\to \mathbb{R}^{n\times m}$ are  locally Lipschitz
functions that  admit $\alpha_4\in \mathcal{K}_\infty$ such that
\[|f(x,t,p(t/\alpha))|+|g(x,t,p(t/\alpha))|\; \le\;
\alpha_4(|x|)\] everywhere, and where $p:\mathbb{R}\to \mathbb{R}^d$
for some $d$ is bounded with a bounded first derivative.  The
control functions (i.e. inputs) for (\ref{controlled}) comprise the
set $\mathcal{U}$ of all measurable essentially bounded functions
$\mathbf{u}:[0,\infty)\to \mathbb{R}^m$ with the essential supremum
norm $|\cdot|_\infty$.  We assume throughout this section that
Assumption \ref{as1} holds for some function $V\in C^1$ and that
$\alpha_1\in \mathcal{K}_\infty$ and $c_a$ from Assumption \ref{as1}
are also such that\smallskip
\begin{itemize}\addtolength{\itemsep}{-0.2\baselineskip}
\item[$A_5$]$|V_x(x,t,p(t/\alpha))|\le c_a\sqrt{\alpha_1(|x|)}$

\item[$A_6$] $|g(x,t,p(t/\alpha))|\le
c_a\left(1+\sqrt[4]{\alpha_1(|x|)}\right)$
\end{itemize}\smallskip
hold for all $t\ge 0$, $\alpha>0$, and $x\in \mathbb{R}^n$. Notice
that $A_5$ reduces to a linear growth condition when
$\alpha_1(x)=|x|^2$ and so  automatically holds in the classical
case where $V$ has the form $x^\top P(t)x$ for a suitable bounded
positive definite matrix. We  show that when Assumption \ref{as1}
and  $A_5$-$A_6$ hold, and when the constant $\alpha>0$ is
sufficiently large, the control system (\ref{controlled}) satisfies
the input-to-state stable (ISS) property and admits the ISS Lyapunov
function (\ref{5}).  We first recall the relevant ISS definitions
from {}\cite{ELW00, MM05, S89, S00}.

For a general locally Lipschitz
control affine system
\begin{equation}
\label{genc} \dot x=h(x,t)+J(x,t)u
\end{equation}
where $h$ and $J$ are uniformly state bounded (as defined in Section
\ref{defs}), and for given values $t_o\ge 0$, $x_o\in \mathbb{R}^n$,
and $\mathbf{u}\in \mathcal{U}$, we let $t\mapsto \phi(t; t_o, x_o,
\mathbf{u})$ denote the unique maximal solution of the initial value
problem \[\dot x(t)=h(x(t),t)+J(x(t),t)\mathbf{u}(t)\; \; {\rm a.e.
}\; \; t\; ,\; \; \; x(t_o)=x_0\; .\] We always assume all
trajectories $\phi(\cdot; t_o, x_o, \mathbf{u})$ so defined are
defined on all of $[t_o,\infty)$.  Later we specialize to the
controlled system (\ref{controlled}) for fixed constants $\alpha>0$.
\begin{definition}
We say that (\ref{genc}) is {\em input-to-state stable (ISS)}
provided there exist $\beta\in \mathcal{KL}$ and $\gamma\in
\mathcal{K}_\infty$ such that \[\tag{ISS} |\phi(t; t_o, x_o,
\mathbf{u})|\; \; \le\; \;
\beta(|x_o|,t-t_o)+\gamma(|\mathbf{u}|_\infty)\] holds for all $t\ge
t_o$, $t_o\ge 0$, $x_o\in \mathbb{R}^n$, and $\mathbf{u}\in
\mathcal{U}$.\end{definition}

\begin{definition}
A smooth function $W:\mathbb{R}^n\times [0,\infty)\to [0,\infty)$ is
called an {\em ISS Lyapunov function} for (\ref{genc}) provided
there exist functions $\alpha_1,\alpha_2,\chi\in\mathcal{K}_\infty$
and a positive definite function $\alpha_3$ such that
\begin{itemize}
\item[(1)] $\alpha_1(|x|)\le W(x,t)\le \alpha_2(|x|)$\ \ and
\item[(2)]
$|u|\le \chi(|x|)$ implies $W_t(x,t)+W_x(x,t)[h(x,t)+J(x,t)u]\le
-\alpha_3(|x|)$.\end{itemize} hold for all $x\in\mathbb{R}^n$, $t\ge
0$, and $u\in \mathbb{R}^m$.
\end{definition}

The following lemma comes from \cite{ELW00}:

\begin{lemma}
\label{edwards} If (\ref{genc}) has an ISS Lyapunov function, then
it is ISS.\end{lemma}

We prove the following analogue of \cite[Theorem 4]{MMdQ06}:

\begin{theorem}
\label{th3} Assume (\ref{1}) satisfies Assumption \ref{as1}.  Assume
$A_5$-$A_6$ everywhere hold where $c_a$, $V$, and $\alpha_1$ are
chosen as in Assumption \ref{as1}.  Then for each constant
$\alpha>4Tc_a\bar p/c_b$, the dynamics (\ref{controlled}) are ISS
and
\begin{equation}
\label{500} \! \! \! V^\sharp_\alpha(t,x) :=
e^{\frac{\alpha}{T}\displaystyle\int_{\frac{t}{\alpha} -
T}^{\frac{t}{\alpha}} \displaystyle\int_{s}^{\frac{t}{\alpha}}
q(p(l)) dl\,  ds} V(x,t,p(t/\alpha))
\end{equation}
is an ISS Lyapunov function for (\ref{controlled}).
\end{theorem}
\noindent{\em Proof:}\ \  We indicate the changes needed in the
proof of Theorem \ref{th1}. Consider the function $\chi\in
\mathcal{K}_\infty$ defined by
\[
\chi(s):=\frac{c_b\sqrt{\alpha_1(s)}}{2Tc^2_a\left(1+\sqrt[4]{\alpha_1(s)}\right)},
\]
where $\alpha_1$ and $c_a$ are as in Assumption \ref{as1}.  This
function is of class $\mathcal{K}_\infty$ since $s\mapsto
\alpha^{1/4}_1(s)$ and $s\mapsto s^2/(1+s)$ are both
$\mathcal{K}_\infty$.  Our assumptions imply that if $|u|\le
\chi(|x|)$, then
\[
\begin{array}{rcl}
|V_x(x,t,p(t/\alpha))g(x,t,p(t/\alpha)u|&\le &
\dfrac{c_b\alpha_1(|x|)}{2T}\\[1em]
&\le & \dfrac{c_b}{2T}V(x,t,p(t/\alpha))
\end{array}
\]
everywhere.  Define $\hat V$ as in (\ref{v1def}) and $E(t,\alpha)$
as in (\ref{eta}). Then, along any trajectory $x=\phi(t)$ of
(\ref{controlled}) with inputs $\mathbf{u}$ satisfying
$|\mathbf{u}|_\infty\le \chi(|\phi(t)|)$ everywhere, we get
\[ \dot{\hat V} = - q(p(t/\alpha))
\hat V(t,x) + \left[\frac{c_a\bar
p}{\alpha}+\frac{c_b}{2T}\right]\hat V(x,t)\] everywhere and
therefore (by reasoning exactly as in Section \ref{mainp}) also
\[
\begin{array}{rcl}
 \! \! \dot{V}^\sharp_\alpha &\le & \! \! E(t,\alpha)\left[\dfrac{c_a\bar
 p}{\alpha}+\dfrac{c_b}{2T}-\dfrac{1}{T}\displaystyle\int_{t/\alpha-T}^{t/\alpha}q(p(l))dl\right]\hat
 V(x,t)\! \\[1em]
 &\le &
 \! \! E(t,\alpha)\left[\dfrac{c_a\bar
 p}{\alpha}-\dfrac{c_b}{2T}\right]\hat V(x,t)\; \le\;   -\dfrac{c_bE(t,\alpha)}{4T}\hat V(x,t)
\end{array}
\]
when $\alpha>4Tc_a\bar p/c_b$, where the second inequality follows
from $A_4$, and the last inequality follows from our choice of
$\alpha$. We then argue exactly as before to show that
$V^\sharp_\alpha$ is an ISS Lyapunov function for (\ref{controlled})
for  all constants $\alpha>4Tc_a\bar p/c_b$.
 The theorem now follows from Lemma
\ref{edwards}. \hfill$\blacksquare$

\begin{remark} \label{lastrm}Theorem \ref{th3} readily applies to our friction
example from Section \ref{revisit} above,  showing that  (\ref{fl})
 is actually an ISS Lyapunov function for
the slowly time-varying controlled friction dynamic
\begin{equation} \label{cfrictionexample}
\begin{array}{rcl}
\! \! \! \dot{x}_{1} &=&x_{2}\\ \dot{x}_{2} &=&-\sigma
_{1}(t/\alpha)x_{2}-k(t)x_{1}+g(x,t,t/\alpha)u -\left\{ \sigma
_{2}(t/\alpha)+\sigma _{3}(t/\alpha)e^{-\beta _{1}\mu
(x_{2})}\right\} {\rm sat} (x_{2})\end{array}\end{equation} for any
$g$ satisfying our assumption $A_6$ for some $c_a>0$ and with
$\alpha_1(s)=s^2/2$, provided the constants $c_a>0$ and $\alpha>0$
are sufficiently large.  Similar extensions can be made for the
other  examples we considered above.
\end{remark}
\section{Conclusion}
\label{concl}


We provided general conditions under which slowly time-varying
systems are uniformly globally asymptotically stable and
input-to-state stable with respect to general perturbations, thus
extending \cite{PA02b, S94} to situations where the corresponding
frozen dynamics are not necessarily exponentially stable.  Moreover,
we provided new methods for constructing explicit closed form strict
ISS Lyapunov functions for slowly time-varying control systems in
terms of a suitable family of generalized Lyapunov like functions
for the frozen dynamics.  This is significant because  Lyapunov
functions play essential roles in robustness analysis and controller
design.

We conjecture that our work can be extended to  systems that are
subjected to actuator or measurement errors, or which  are
components of larger controlled hybrid dynamical systems.  It would
also be of interest to extend our work to slowly time-varying
systems with outputs and  to construct corresponding input-to-output
stable (IOS) Lyapunov functions; see \cite{SW99, SW01} for further
background on systems with outputs and \cite{MM05a} for some first
results on constructing IOS Lyapunov functions for non-autonomous
systems in terms of given nonstrict Lyapunov functions.
We leave these extensions for future papers.

\smallskip

\section*{Acknowledgements}
This research was supported by NSF Grant 0424011. This work was done
while the first author visited Louisiana State University (LSU). He
thanks LSU for the kind hospitality he enjoyed during this period.
Both authors thank Marcio de Queiroz and Patrick De Leenheer for
illuminating discussions.



\begin{thebibliography}{999999}

\bibitem[A]{A99}
D. Angeli,  Input-to-State Stability of PD-controlled robotic
systems, {\em Automatica 35} (1999), 1285-1290.


\bibitem[AS]{AS99} D. Angeli and E.D. Sontag, Forward
completeness, unboundedness observability, and their Lyapunov
characterizations, {\em Systems and Control Letters 38} (1999),
 209-217.

\bibitem[ASW]{ASW00} D. Angeli, E.D. Sontag, and Y. Wang, A
characterization of integral input to state stability, {\em IEEE
Trans.  Automatic  Control  45} (2000), 1082-1097.




\bibitem[DDNZ]{DDNZ00} M.S. de Queiroz,  D.M. Dawson, S. Nagarkatti, and
F. Zhang, {\it Lyapunov-Based Control of Mechanical Systems},
Birkh\"{a}user, Cambridge, MA, 2000.

\bibitem[ELW]{ELW00}
H. Edwards,  Y. Lin, and Y. Wang,  On input-to-state stability for
time-varying nonlinear systems, in {\em Proceedings of the  39th
IEEE Conf. on Decision and Control}, Sydney, Australia, December
2000,  3501-3506.





\bibitem[K]{K02}
H. Khalil, {\it Nonlinear Systems, Third Edition}, Prentice Hall,
Englewood Cliffs, NJ,  2002.






\bibitem[MM1]{MM05a} M. Malisoff  and F. Mazenc, Further constructions of strict
Lyapunov functions for time-varying systems, in {\em Proceedings of
the American Control Conference},  Portland, OR, June, 2005,
1889-1894.

\bibitem[MM2]{MM05}
M.  Malisoff and F. Mazenc, Further remarks on strict input-to-state
stable Lyapunov functions for time-varying systems, {\em Automatica
41} (2005),  1973-1978.





\bibitem[M]{M03}
F.  Mazenc, Strict Lyapunov functions for time-varying systems, {\em
Automatica  39} (2003),  349-353.

\bibitem[MM3]{MM06} F. Mazenc and M. Malisoff, Further
constructions of control-Lyapunov functions and stabilizing
feedbacks for systems satisfying the Jurdjevic-Quinn conditions,
 {\em IEEE Trans.  Automatic  Control 51} (2006), 360-365.

\bibitem[MMD]{MMdQ06}F. Mazenc,  M. Malisoff, and M. de Queiroz, Further
results on strict Lyapunov functions for rapidly time-varying
nonlinear systems, {\em Automatica 42} (2006), 1663-1671.














  \bibitem[PA1]{PA02a}
J. Peuteman  and D. Aeyels, Exponential stability of nonlinear
time-varying differential equations and partial averaging, {\it
Mathematics of Control, Signals, and Systems 15} (2002),  42-70.

\bibitem[PA2]{PA02b}
J.  Peuteman  and D. Aeyels, Exponential stability of slowly
time-varying nonlinear systems, {\it Mathematics of Control,
Signals, and Systems 15} (2002),  202-228.



\bibitem[S1]{S94}
V. Solo, On the stability of slowly time-varying linear systems,
{\it Mathematics of Control, Signals, and Systems 7} (1994),
331-350.



\bibitem[S2]{S89}
E.D. Sontag,   Smooth stabilization implies coprime factorization,
{\em IEEE Trans.   Automatic    Control 34} (1989), 435--443.


\bibitem[S3]{S00}
E.D. Sontag, The ISS philosophy as a unifying framework for
stability-like behavior, in {\em Nonlinear Control in the Year 2000,
Volume 2}, A. Isidori, F. Lamnabhi-Lagarrigue, and W. Respondek
(eds.), Springer-Verlag, Berlin, 2000, 443-468.



\bibitem[SW1]{SW99} E.D. Sontag and Y. Wang, Notions of input to
output stability, {\em Systems and Control Letters 38} (1999),
235--248.

\bibitem[SW2]{SW01} E.D. Sontag and Y. Wang, Lyapunov
characterizations of input to output stability, {\em SIAM J. Control
and Optimization 39} (2001), 226-249.




\end{thebibliography}
\end{document}